\documentclass[11pt,a4paper]{article}
\usepackage{ifthen,latexsym,amssymb,amsmath,bbm,fixmath}
\usepackage[nobysame,initials]{amsrefs}
\usepackage[shortlabels]{enumitem} 
\usepackage{hyperref}
\usepackage{xcolor}

\newcommand{\C}[1]{{\protect\mathcal{#1}}}

\newcommand{\I}[1]{{\mathbbm #1}}
\renewcommand{\O}[1]{\overline{#1}}

\newcommand{\ceil}[1]{\lceil #1\rceil}
\newcommand{\e}{\varepsilon}
\newcommand{\floor}[1]{\lfloor #1\rfloor}

\renewcommand{\mid}{:}

\renewcommand{\ge}{\geqslant}
\renewcommand{\le}{\leqslant}

\newif\ifnotesw\noteswtrue

\newcommand{\hide}[1]{}

\newcommand{\PST}[1]{\ifthenelse{\equal{#1}{}}
{\cite{PikhurkoSliacanTyros19}}
{\cite{PikhurkoSliacanTyros19}*{#1}}}

\newcommand{\beq}[1]{\begin{equation}\label{#1}}
\newcommand{\eeq}{\end{equation}}

\newtheorem{theorem}{Theorem}
\newtheorem{lemma}[theorem]{Lemma}

\newtheorem{corollary}[theorem]{Corollary}

\newtheorem{proposition}[theorem]{Proposition}

\newcommand{\bpf}[1][Proof.]{\smallskip\noindent{\it #1} }
\newcommand{\qed}{\nolinebreak\mbox{\hspace{5 true pt}%
  \rule[-0.85 true pt]{3.9 true pt}{8.1 true pt}}}
\newcommand{\epf}{\qed \medskip}

\newcommand{\K}[2]{\kappa(#2;#1)}

\begin{document}


\newcommand{\OK}{\O K}

\newcommand{\R}{R_{3,5}}

\newcommand{\ER}[2]{\mathrm{ER}_{#1}(#2)}

\newcommand{\er}[2]{\mathrm{er}_{#1}(#2)}

\newcommand{\BF}[1]{\mathcal{B}_{#1}} 

\newcommand{\BB}[2]{\mathrm{BB}_{#1}(#2)}

\newcommand{\bb}[2]{\er{#1}{\BF{#2}}}

\newcommand{\Forb}[1]{\mathrm{Forb}_{#1}}

\newcommand{\blow}[2]{#1(\!(#2)\!)}

\newcommand{\dedit}{d_{\mathrm{edit}}}

\newcommand{\m}{m} 


\newcommand{\OO}[1]{\overline{#1}}
\newcommand{\pp}{P}
\newcommand{\barL}{R_{3,3,3}}

\newcommand{\PV}[1]{\ifthenelse{\equal{#1}{}}
{\cite{PikhurkoVaughan13}}
{\cite{PikhurkoVaughan13}*{#1}}}

\title{On problems of Erd\H os and Baumann--Briggs on minimising the density of $s$-cliques in graphs with forbidden subgraphs}
\author{Levente Bodn\'ar and Oleg Pikhurko\\Mathematics Institute and DIMAP\\
University of Warwick\\
Coventry CV4 7AL, UK}
\maketitle

\begin{abstract}
Using flag algebras, we prove that the minimum density of $8$-cliques in a large graph without an independent set of size $3$ is
$\frac{491411}{268435456}+o(1)$, thus resolving a new case of an old problem of Erd\H os [\emph{Magyar Tud.\ Akad.\ Mat.\ Kutat\'o Int.\ K\"ozl.} 7 (1962) 459--464]. Also, we establish some other results of this type; for example,  we show that the minimum $s$-clique density in a large graph with no independent set of size 3 nor an induced 5-cycle is $2^{1-s}+o(1)$ when $s=4,5,6$. For each of these results, we also describe the structure of all extremal and almost extremal graphs of large order~$n$.

These results are applied to give an asymptotic solution to a number of cases of the problem of Baumann and Briggs [\emph{Electron. J. Combin.} 32 (2025) P1.22] which asks for the minimum number of $s$-cliques in an $n$-vertex graph in which every $k$-set spans a $t$-clique.
\end{abstract}

\section{Introduction}
\label{se:Intro}

As usual, a \emph{graph} $G$ is a pair
$(V(G),E(G))$, where $V(G)$ is the \emph{vertex set} and the \emph{edge
set} $E(G)$ consists of unordered pairs of vertices. 
We let $\OO G:=(V(G),\binom{V(G)}{2}\setminus E(G))$
denote its \emph{complement} and $v(G):=|V(G)|$
denote its \emph{order}. For graphs $F$ and $G$,
let $\pp(F,G)$ be the number of $v(F)$-subsets of $V(G)$
that induce in $G$ a subgraph isomorphic to $F$.%
\hide{; further, for $v(G)\ge v(F)$, define the 
\emph{density of $F$ in $G$} to be
 \beq{density}
 p(F,G):=\pp(F,G)\,{v(G)\choose v(F)}^{-1}.
 \eeq}
 Let $K_\ell$ denote the \emph{$\ell$-clique}, that is, 
 the complete graph with $\ell$
 vertices. Thus, for example,  $\OK_\ell$ is the empty graph on $\ell$ vertices. 

Given a family $\C F$ of graphs and an integer $s\ge 2$, we consider the following extremal function:
 \beq{eq:ErGeneral} 
 \ER{s}{n,\C F}:=\min\{\pp(K_s,G)\mid v(G)=n,\ \forall\, F\in\C F\ \pp(F,G)=0\},\quad \mbox{for $n\in\I N$},
 \eeq
 the minimum number of $s$-cliques in an $n$-vertex graph $G$ which is \emph{$\C F$-free} (that is, no $F\in\C F$ can appear in $G$ as an induced subgraph). If the graph family $\C F=\{F\}$ consists of a single graph, then we write $\ER{s}{n,F}$ for $\ER{s}{n,\{F\}}$, etc.
 
The case $s=2$ of~\eqref{eq:ErGeneral} includes, after passing to graph complements,  the classical Tur\'an function. This was fully resolved by Tur\'an~\cite{Turan41} when any given clique is forbidden (see also Mantel~\cite{Mantel07}), while the Erd\H os--Stone Theorem~\cite{ErdosStone46} determines it within an additive $o(n^2)$ error term for any forbidden family (with the corresponding stability result proved by Erd\H os~\cite{Erdos67a} and Simonovits~\cite{Simonovits68}).
Here we are interested in the case $s\ge 3$.  

The special case of~\eqref{eq:ErGeneral} when $s=3$ and $\C F=\{\OK_3\}$ 
was resolved by Goodman~\cite{Goodman59} (for all even $n$), with  Lorden~\cite{Lorden62} determining the exact value for all odd $n$. Erd\H os~\cite{Erdos62a} asked for the value of $\ER{s}{n,\OK_\ell}$ and made a too optimistic conjecture about it which was disproved by Nikiforov~\cite{Nikiforov01}.
It is easy to show that the limit
  \beq{ckl}
  \er{s}{\C F}:=\lim_{n\to\infty} \frac{\ER{s}{n,\C F}}{{\binom{n}{s}}}
  \eeq
  exists for every graph family $\C F$ (see e.g.~\PST{Lemma~2.2} for a proof). However, the value of $\er{s}{\OK_{\ell}}$ for any other non-trivial pair $(s,\ell)$ apart from those mentioned above was not known for a long while. 

After the introduction of the flag algebra method by Razborov~\cites{Razborov07,Razborov10}, a number of new cases have been solved. 
In order to state these results, we need some further definitions. For an integer $m\ge 0$, we denote $[m]:=\{0,\dots,m-1\}$.
For $m\ge 3$, let $C_m$ denote the \emph{$m$-cycle} whose vertex set is $[m]$ with the edge set consisting of the pairs of consecutive vertices modulo~$m$.
Let $\R$ denote the graph on $[13]=\{0,\dots,12\}$ in which $x,y\in V(\R)$ are adjacent if $x-y$ modulo 13 is in $\{1,5,8,12\}$. It is easy to check that $\R$ does not contain $K_3$ nor $\OK_5$ as induced subgraphs. (In fact, $\R$ is the unique $(3,5)$-Ramsey graph on $13$ vertices.) Let $\barL$ be the graph whose vertex set consists of all even-sized subsets of $[5]$, where $X,Y\in V(\barL)$ are adjacent if $|X\bigtriangleup Y|=2$, that is, their symmetric difference has 2 elements. Thus $\barL$ is a 10-regular graph with 16 vertices. It does not contain $K_6$ nor $\OK_3$. In fact, $\barL$ is the complement of the Clebsch graph (which appears in extremal $(3,3,3)$-Ramsey colourings). 
Given a graph $H$, which by relabelling its vertices is assumed to have vertex set
$[m]$,
and disjoint
sets $V_0,\dots,V_{m-1}$,
let the \emph{expansion} 
$\blow{H}{V_0,\dots,V_{m-1}}$ be the graph on $V_0\cup\dots\cup V_{m-1}$ obtained
by
putting the complete graph on each $V_i$ and
putting, for each edge $\{i,j\}\in E(H)$, the complete bipartite graph between
$V_i$ and~$V_j$. 
An expansion is \emph{uniform} if  $\big|\,|V_i|-|V_j|\,\big|\le 1$ for any $i,j\in [m]$.  
If we consider expansions in terms of complements, then it amounts
to blowing up each vertex $i$ of $\OO H$ by a factor of $|V_i|$. 
The \emph{edit distance} 
between two graphs of the same order is the smallest number of adjacencies that have to be changed in one graph to make it isomorphic to the other. We call a graph \emph{$\ER{s}{n,\C F}$-extremal} if it has $n$ vertices, is $\C F$-free and has $\ER{s}{n,\C F}$ $s$-cliques. We call a graph $G$ of order $n$ (or rather a sequence of graphs with $n\to\infty$) \emph{almost $\ER{s}{n,\C F}$-extremal} if $G$ is $\C F$-free and $\pp(K_s,G)=(\er{s}{\C F}+o(1))\binom ns$. If $s$, $\C F$ and $n$ are understood, we may just say \emph{(almost) extremal}.

Now we are ready to state the following theorem that collects some of the results established by Das, Huang, Ma, 
Naves, and Sudakov~\cite{DasHuangMaNavesSudakov13}, Parczyk, Pokutta, Spiegel and Szab\'o~\cite{ParczykPokuttaSpiegelSzabo23}, and Pikhurko and Vaughan~\PV{} (in addition to the results from \cites{Goodman59,Lorden62} for which the stated stability property can be easily established by modern methods).

\begin{theorem}[\cites{DasHuangMaNavesSudakov13,Goodman59,Lorden62,ParczykPokuttaSpiegelSzabo23,PikhurkoVaughan13}]
\label{th:ER}
Let $(s,\ell)$ be one of the pairs of integers listed below and let $H$ be the following graph:
\begin{enumerate}[(a)]
\item\label{it:K} $\OK_{\ell-1}$ if $s=3$ and $3\le \ell\le 7$, 
\item\label{it:C} $C_5$ if $s\in\{4,5\}$ and $\ell=3$, 
\item\label{it:L} $\barL$ if $s\in\{6,7\}$ and $\ell=3$,
\item\label{it:R} $\R$ if $s=4$ and $\ell=5$. 
\end{enumerate}
Then every almost $\ER{s}{n,\OK_\ell}$-extremal graph is $o(n^2)$-close in the edit distance to a uniform expansion of~$H$. Moreover, 
for all sufficiently large $n$,
every $\ER{s}{n,\OK_\ell}$-extremal graph contains an expansion of $H$ as a spanning subgraph. 
\end{theorem}

Formally, the first conclusion of Theorem~\ref{th:ER} states that for every $\e>0$ there are $\delta>0$ and $n_0$ such that every $\OK_\ell$-free graph $G$ with $n\ge n_0$ vertices and $\pp(K_s,G)\le (1+\delta)\ER{s}{n,\OK_\ell}$ is within edit distance  $\e{\binom{n}{2}}$ from a uniform expansion of $H$. When the meaning is clear, we will be using the shorter $o(1)$-version and refer to this property (that all almost extremal $n$-vertex graphs are within edit distance $o(n^2)$ from each other) as \emph{stability}. In Cases~\ref{it:C}--\ref{it:R}, it holds that, in fact, every large $\ER{s}{n,\OK_\ell}$-extremal graph is an expansion of~$H$ (without any extra edges added). 

Our first result expands the list of the solved cases of the Erd\H os problem.

\begin{theorem}\label{th:NewER} Every almost $\ER{8}{n,\OK_3}$-extremal graph is $o(n^2)$-close in the edit distance to a uniform expansion of~$\barL$. Moreover, for every sufficiently large $n$,
every $\ER{8}{n,\OK_3}$-extremal graph is an expansion of~$\barL$. 
\end{theorem}

A routine application of the Inclusion-Exclusion Principle gives that the limit density of $K_s$ in growing uniform expansions of $\barL$ is
 \beq{eq:KsL}
  \frac{5\cdot 5^{s-1} + 10\cdot 4^{s-1} - 30\cdot 3^{s-1} + 20\cdot 2^{s-1}-4}{16^{s-1}}.
 \eeq
 Thus, it follows from Theorem~\ref{th:NewER} that
 \[
 \er{8}{\OK_3}=\frac{491411}{268435456}.
 \]

A related extremal problem was posed by Baumann and Briggs~\cite{BaumannBriggs25}, who asked for the value of
$
\ER{s}{n,\BF{k,t}},
$
 where \beq{eq:Fkt}
 \BF{k,t}:=\{F\mid v(F)=k,\ \pp(K_t,F)=0\}
 \eeq
  consists of all $k$-vertex graphs containing no $t$-clique. In other words, $\ER{s}{n,\BF{k,t}}$ is the minimum number of $s$-cliques in an $n$-vertex graph such that every $k$-set of vertices spans at least one $t$-clique.
 Here, the case $s=2$ falls into the remit of the Erd\H os--Stone Theorem; 
 in fact, the $\ER{2}{n,\BF{k,t}}$-problem  was fully resolved in~\cites{BaumannBriggs25,HoffmanJohnsonMcDonald15,NobleJohnsonHoffmanMcDonald17}. Also, if $t=2$ then we forbid only $\OK_k$ and thus we get the function $\ER{s}{n,\OK_k}$ of Erd\H os; this was one of the motivations in~\cite{BaumannBriggs25} for introducing the function~$
\ER{s}{n,\BF{k,t}}
$. 
 
In this paper, we relate the above two problems as follows. Let $\BF{k,t}'$ consist of all graphs that admit a $k$-vertex expansion not containing a $t$-clique. Clearly, $\BF{k,t}'\supseteq \BF{k,t}$. Also, since we allow empty parts when taking expansions, if an induced subgraph of a graph $H$ is in $\BF{k,t}'$ then $H$ itself is also in~$\BF{k,t}'$.
Our translation of results on the Erd\H os function to ones about the Baumann--Briggs function goes by taking $\C F=\{\OK_\ell\}$ in the following proposition.

\begin{proposition}\label{pr:Reduction2} Let $k\ge t\ge 2$ be integers and 
$\C F$ be a subfamily of $\BF{k,t}'$. If an integer $s\ge 2$ and a graph $H\not\in\BF{k,t}'$ satisfy that some (not necessarily uniform) expansions of $H$ give almost extremal graphs for the $\er{s}{\C F}$-problem then
 \beq{eq:Reduction2=}
 \er{s}{\BF{k,t}}=\er{s}{\C F}.
\eeq
 Moreover, as $n\to\infty$, every almost $\ER{s}{n,\BF{k,t}}$-extremal graph is $o(n^2)$-close in the edit distance to some almost extremal $\ER{s}{n,\C F}$-graph.
 \end{proposition}

In order to apply these results to the Baumann--Briggs problem, we have to determine, for each of the graphs $\OK_{\ell}$, $C_5$, $\barL$ and $\R$, the set of pairs $(k,t)$ for which the graph is in $\BF{k,t}'$. By the monotonicity in $k$, it is enough to know the following function:
for a graph $H$ and an integer $t\ge 2$, define $\K{H}{t}$ to be the maximum $k$ such that there is a $k$-vertex expansion of $H$ without a $t$-clique (that is, $H\in\BF{k,t}'$). Thus the assumptions that $H\not\in \BF{k,t}'$ and $\C F\subseteq \BF{k,t}'$ of Proposition~\ref{pr:Reduction2} can be restated as
 \begin{equation}\label{eq:Reduction2}
 \K{H}{t}< k\le \min\{\K{F}{t}\mid F\in\C F\}.
 \end{equation}
 As an illustration, let us determine this function for edgeless graphs, postponing  the analogous results
for the remaining graphs to Section~\ref{se:BB}.

\begin{lemma}\label{lm:OK} For integers $t\ge 2$ and $\m\ge 1$, we have $\K{\OK_{\m}}{t}=\m(t-1)$.\end{lemma}
\bpf A uniform expansion of $\OK_{\m}$ with each part having $t-1$ vertices has no $K_t$; this shows that $\K{\OK_{\m}}{t}\ge \m(t-1)$. On the other hand, if $k>(t-1)\m$ then every $k$-vertex expansion of $\OK_{\m}$ has a part with at least $t$ vertices by the Pigeonhole Principle, giving a $t$-clique.\epf

By combining this (for $m=\ell-1$ and $m=\ell$) with Proposition~\ref{pr:Reduction2} (for $\C F=\{\OK_\ell\}$ and $H=\OK_{\ell-1}$) and Theorem~\ref{th:ER}\ref{it:K} whose conclusion is (trivially) true  also if $\ell=2$, 
we obtain by~\eqref{eq:Reduction2} the following result.

\begin{corollary}\label{cr:1} If integers $k\ge t\ge 2$ and $2\le \ell\le 7$ satisfy $(t-1)(\ell-1)<k\le (t-1) \ell$ then
 \beq{eq:Kl}
 \er{3}{\BF{k,t}}=\frac{1}{(\ell-1)^2},
  \eeq
  and every almost $\ER{3}{n,\BF{k,t}}$-extremal graph is $o(n^2)$-close in the edit distance to a uniform expansion of~$\OK_{\ell-1}$.\qed
  \end{corollary}
  
Note the following easy result valid for all integers $s\ge 2$ and $k\ge t\ge 2$ (which for $s=3$ follows from Corollary~\ref{cr:1} with $\ell=2$):
\beq{eq:1}
\er{s}{\BF{k,t}}=1,\quad\mbox{if $k\le 2t-2$}.
\eeq
Indeed, if $k\le 2t-2$ then the complement of any $n$-vertex $\BF{k,t}$-free graph $G$ does not contain the complete bipartite graph with parts of sizes $t-1$ and $k-t+1$, and thus has at most $O(n^{2-1/(k-t+1)})=o(n^2)$ edges as $n\to\infty$ by the well-known bound of K{\H o}v{\'a}ri, S{\'o}s and Tur{\'a}n~\cite{KovariSosTuran54}. Consequently, at most $|E(\OO G)|\binom{n-2}{s-2}=o(n^s)$ $s$-sets of $G$ fail to span a clique, proving~\eqref{eq:1}.
  
Also, Corollary~\ref{cr:1} 
gives that  $\bb{3}{5,3}=1/4$ and all almost extremal graphs are close to a uniform expansion of $\OK_2$, giving an asymptotic answer to~\cite{BaumannBriggs25}*{Question~19} that asks for the value of $\ER{3}{n,\BF{5,3}}$. See also Theorem~\ref{th:53} that describes the structure of all large extremal graphs for this problem.

There are other cases of the Baumann--Briggs problem which we could asymptotically solve by forbidding one further graph in addition to an independent set. Here we state two such theorems, which we believe to be  of independent interest (with their consequences for the Baumann--Briggs problem discussed in Section~\ref{se:BB}).

\begin{theorem}\label{th:ForbC5} Take any $s\in \{4,5,6\}$. Then
$\er{s}{\{\OK_3,C_5\}}=2^{1-s}$ and every almost extremal  graph with $n\to\infty$ vertices is $o(n^2)$-close in the edit distance to a uniform expansion of $\OK_{2}$. Moreover, for all sufficiently large $n$, every extremal graph contains a uniform expansion of~$\OK_{2}$ as a spanning subgraph and it holds that $\ER{s}{n,\{\OK_3,C_5\}}=\binom{\floor{n/2}}s+\binom{\ceil{n/2}}s$.
\end{theorem}

Let $C_5'$ be the 6-vertex graph which is the disjoint union of a $5$-cycle $C_5$ and an isolated vertex.

\begin{theorem}\label{th:734} It holds that $\er{4}{\{\OK_4,C_5'\}}=1/27$ and every almost $\ER{4}{n,\{\OK_4,C_5'\}}$-extremal graph can be transformed into a uniform expansion of $\OK_3$ by changing at most $o(n^2)$ adjacencies. Moreover, for all sufficiently large $n$, every $\ER{4}{n,\{\OK_4,C_5'\}}$-extremal graph has a uniform expansion of~$\OK_3$ as a spanning subgraph and it holds that $\ER{4}{n,\{\OK_4,C_5'\}}=\sum_{i=0}^2 \binom{\floor{(n+i)/3}}4$.
\end{theorem}

\section{Consequences for  the Baumann--Briggs problem}
\label{se:BB}

In this section, we analyse when each of the graphs $C_5$, $\R$ and $\barL$ is in $\BF{k,t}'$ and state the corresponding consequences for the Baumann--Briggs problem. The following easy lemma will be used a few times here.

\begin{lemma}\label{lm:General} Let $H$ be an $\ell$-vertex graph  and let $\omega,d\ge 1$ be integers such that every vertex of $H$ is contained in exactly $d$ $\omega$-cliques. If $H\in \BF{k,t}'$ (that is, $H$ admits a $k$-vertex $K_t$-free expansion) for some integers $k\ge t\ge 2$ then $\ell (t-1)\ge \omega k$.\end{lemma}

\bpf Let $\C C$ denote the collection of the vertex sets of  $\omega$-cliques in $H$. Clearly, we have that $|\C C|=d\ell /\omega$.
Take an expansion $\blow{H}{V_0,\dots,V_{\ell-1}}$ of $H$ with $k$ vertices in which  every clique has at most $t-1$ vertices. 
Thus, for every $C\in \C C$, we have that $t-1\ge \sum_{u\in C} |V_u|$. Summing this over all $\omega$-cliques in $H$ we obtain:
\[
\frac{d\ell}{\omega}\, (t-1)=\sum_{C\in \C C} (t-1)\ge \sum_{C\in\C C} \sum_{u\in C} |V_u|
=\sum_{u\in V(H)} |V_u|\cdot |\{C\in\C C\colon C\ni u\}|= kd,
\]
giving the claimed inequality.\epf

Recall that for a graph $H$ and integer $t\ge 2$, we can define $\K{H}{t}$ as the maximum $k$ such that $H\in \BF{k,t}'$. Alternatively, $\K{H}{t}+1$ is the smallest integer $m$ such that every expansion of $H$ with $m$ vertices has a $t$-clique.
For example, $\K{H}{2}$ is the independence number of $H$. For convenience, let us agree that $\K{H}{1}:=0$, which matches what the definition formally  gives in this case.

\begin{lemma}\label{lm:C5} For any integer $t\ge 1$, it holds that $\K{C_5}{t}=\floor{5(t-1)/2}$.
\end{lemma}
\bpf The case $t=1$ follows from the definition, so assume that $t\ge 2$. Let $t=2m+c+1$ for an integer $m\ge 0$ and $c\in \{0,1\}$. Let $b:=0$ if $c=0$ and $b:=2$ if $c=1$. Thus we have to show that $\K{C_5}{2m+c+1}=5m+b$.

For the lower bound, we have to find a $(5m+b)$-vertex expansion of $C_5$ 
that
has no $(2m+c+1)$-clique. We  take $b$ parts of size  $m+1$ and $5-b$ parts of size $m$ with the  restriction that if $c=1$ (when $b=2$) then  the two strictly larger parts (or, more precisely, their indices) are not adjacent in $C_5$. Since $C_5$ is triangle-free, every clique in the expansion intersects at most two parts which then have to be adjacent. If $c=0$ (resp.\ $c=1$) then any two adjacent parts get at most  $2m$ (resp.\ $2m+1$) vertices, as desired.

On the other hand, Lemma~\ref{lm:General} applied to $C_5$ (with $\omega=d=2$), shows that every $C_5$-expansion with $5m+b+1$ vertices has a clique with $2m+c+1$ vertices. Indeed, otherwise it would hold that $5(2m+c)\ge 2(5m+b+1)$, contradicting the definition of $b$.\epf

The combination of Theorem~\ref{th:ER}\ref{it:C} and  Proposition~\ref{pr:Reduction2} applies to those pairs $(k,t)$ such that $C_5\not\in \BF{k,t}'$ and $\OK_3\in \BF{k,t}'$.
Using the values of these functions returned by Lemma~\ref{lm:OK} and Lemma~\ref{lm:C5} we obtain by~\eqref{eq:Reduction2} the following. 


\begin{corollary}\label{cr:C5} If $s\in \{4,5\}$ and integers $k\ge t\ge 2$ satisfy 
$5(t-1)/2< k\le 3(t-1)$
then $\er{s}{\BF{k,t}}=(2^s-1)/5^{s-1}$ and every almost 
$\ER{s}{n,\BF{k,t}}$-extremal graph 
is $o(n^2)$-close in the edit distance to a uniform expansion of $C_5$.\qed\end{corollary}

Next, we determine the function $\K{\R}{t}$ for $\R$, the unique $(3,5)$-Ramsey graph on $13$ vertices (defined in Section~\ref{se:Intro}).

\begin{lemma}\label{lm:R} For every integer $m\ge 0$, it holds that $\K{\R}{2m+1}=13m$ and
$\K{\R}{2m+2}=13m+4$.
\end{lemma}
\bpf For $c\in \{0,1\}$, let $b:=0$ if $c=0$ and $b:=4$ if $c=1$. Thus we have to show that $\K{\R}{2m+c+1}=13m+b$. Ignoring the trivial case $(m,c)=(0,0)$, assume that $2m+c+1\ge 2$.

For the lower bound, it is enough to find a $(13m+b)$-vertex expansion $G=\blow{\R}{V_0,\dots,V_{12}}$ having no clique with $2m+c+1$ vertices. We take $b$ parts of size $m+1$ and $13-b$ parts of size $m$, with the restriction that if $c=1$ (when $b=4$) then the 4 strictly larger parts correspond to an independent set in $\R$ (e.g.\ $\{0,2,4,6\}$). Since $\R$ is triangle-free, for any clique $C$ in $G$ there are two adjacent $i,j$ in $\R$ with $C\subseteq V_i\cup V_j$. We have that $|V_i|+|V_j|$ is at most $2m$ (resp.\ $2m+1$) if $c=0$ (resp.\ $c=1$), as desired.

For the other direction, take any $(13m+b+1)$-vertex expansion $G=\blow{\R}{V_0,\dots,V_{12}}$. If $c=0$ then
Lemma~\ref{lm:General} (with $\omega=2$ and $d=4$) implies that $G$ has  a $(2m+1)$-clique for otherwise $13\cdot 2m\ge 2\cdot (13m+1)$, a contradiction.

The case $c=1$ requires a different argument. Suppose on the contrary that $G$ has no $(2m+2)$-clique. Let $\mu:=\max\{|V_i|-m\colon i\in[13]\}>0$. By the vertex-transitivity of $\R$, assume that $|V_0|=m+\mu$. For each  of the 4 neighbours $i$ of $0$ in $\R$ (namely for $i\in \{1,5,8,12\}$), we have that $|V_i|\le (2m+1)-|V_0|\le m-\mu+1$. The remaining 8 vertices of $\R$ contain a perfect matching $M$, e.g.\ with edges $\{2,10\}$, $\{3,11\}$, $\{ 4,9\}$ and $\{6,7\}$. For each $\{i,j\}\in M$, we have $|V_i|+|V_j|\le 2m+1$ and thus the 8 parts non-adjacent to $V_0$ contain at most $4(2m+1)$ vertices. We conclude that
 \[
 13m+5\le (m+\mu)+4(m-\mu+1)+4(2m+1)=13m-3\mu+8.
 \]
 Thus necessarily $\mu=1$ and we have at least 5 parts of size $m+1$. But some two of these parts are adjacent (as $\R$ has no independent set of size $5$), giving a $(2m+2)$-clique, a contradiction.\epf

Interestingly, if we consider those $k\ge t\ge 3$ for which the combination of Theorem~\ref{th:ER}\ref{it:R} and Proposition~\ref{pr:Reduction2} 
gives that almost $\er{4}{\BF{k,t}}$-extremal graphs are close to expansions of $\R$ then we get the empty set. Indeed, every such pair $(k,t)$ has to satisfy that  $\OK_{5}\in \BF{k,t}'$ (i.e.\ $k\le 5(t-1)$) and $\R\not\in \BF{k,t}'$ (i.e.\ $k> \K{\R}{t}$). For $t=2m+1\ge 3$, we have $\K{\R}{t}-5(t-1)=3m\ge 0$; for $t=2m+2\ge 3$, this difference is $3m-1\ge 0$. Thus there are no such pairs with $t\ge 3$.
\hide{
(* when is K35 BB-extremal *)

tmax[k_] := Module[{m, b, c},
  b = Mod[k - 1, 13] + 1;
  m = (k - b)/13;
  c = If[b <= 4, 1, 2];
  (* Return[{k,m,b,c}] *)
  Return[2 m + c]
  ]

tmin[k_] := Ceiling[k/5]+1

Table[{tmax[k] - tmin[k]}, {k, 2, 100}]
}%

Let us turn our attention to $\barL$, the complement of the Clebsch graph.

\begin{lemma}\label{lm:L} Define $b=b(c)$ to be $0$, $2$, $5$, $9$ and $12$ when $c$ is respectively  $0$, $1$, $2$, $3$ and $4$. Then for any integer $m\ge 0$, it holds that $\K{\barL}{5m+c+1}=16m+b$ except in the case $(m,c)=(0,3)$, when $\K{\barL}{4}=8$.
\end{lemma} 
\bpf  Recall that the vertex set of $\barL$  consists of all even-sized subsets of~$[5]$, with two sets being adjacent if their symmetric difference has exactly 2 elements. Note that this graph has many symmetries (we can permute the ground set $[5]$ or take the symmetric difference with any fixed even-sized set $X\subseteq [5]$), which greatly helps in verifying some of the claims below.
 
It will be useful to describe all 5-cliques in $\barL$. 
For $i\in[5]$, let $C_i\subseteq V(\barL)$ consist of the empty set $\emptyset$ and the four 2-sets $\{i,j\}$ with $j\in [5]\setminus\{i\}$. For distinct $i,j\in [5]$, let $C_{i,j}\subseteq V(\barL)$ consist of the sets $[5]\setminus \{i\}$, $[5]\setminus\{j\}$ and all three 2-subsets of $[5]\setminus\{i,j\}$. Finally, let $C\subseteq V(\barL)$ consist of all sets of size $4$. It is easy to see using the symmetries of $\barL$ that each of the above 16 sets spans a clique in $\barL$ and these are all the 5-cliques. Note that $\barL$ has inclusion-maximal cliques with 4 vertices, so we will have to take these into account too.

Except when $(m,c)=(0,3)$, for the lower bound we have to find a $(16m+b)$-vertex expansion without a clique of order $5m+c+1$. 

For $c=0$ (when $b=0$), we let each part have $m$ vertices. For $c=1$ (when $b=2$), each part has $m$ vertices except for two non-adjacent parts of size $m+1$. For $c=2$ (when $b=5$), each part has $m$ vertices except for five parts of size $m+1$ that correspond to an induced $5$-cycle in~$\barL$. In these three cases, a clique meeting $r\le 5$ parts has at most $rm+c\le 5m+c$ vertices. 

If $c=3$ (when $b=9$) and $m\ge 1$ then, interestingly, we have a part with fewer than $m$ elements. 
Namely, let $|V_X|$ be $m$ if $|X|=4$, $m+1$ if $|X|=2$ and $m-1$ if $X=\emptyset$. For each of the cliques $C$, $C_i$ and $C_{i,j}$ in $\barL$, the total number of vertices in the corresponding parts is respectively at most $5m$, $4(m+1)+(m-1)$ and $2m+3(m+1)$, which is at most the desired upper bound $5m+3$ in each case. A clique meeting at most four parts has at most $4(m+1)\le 5m+3$ vertices, since $m\ge 1$.

It remains to provide a construction for $c=4$ (when $b=12$). 
Here, each part has $m+1$ vertices except each of the four parts indexed by the sets in
$\C A:=\{\emptyset,\{3,4\},[5]\setminus\{4\}, [5]\setminus\{3\}\}$ has $m$ vertices. 
First, let us check that  every $5$-clique in $\barL$ intersects $\C A$. 
 Indeed, the clique $C$ contains e.g.\ $[4]\in\C A$ and each $C_i$ contains $\emptyset\in\C A$. 
 Also, if $\{i,j\}\cap \{3,4\}\ni a$ is non-empty then  $C_{i,j}$ contains $[5]\setminus\{a\}\in\C A$; otherwise $C_{i,j}$ contains $\{3,4\}\in\C A$. Thus a clique meeting five parts has at most $4(m+1)+m$ vertices, while one meeting at most four parts has at most $4(m+1)\le 5m+4$ vertices, as desired.

For the exceptional case $m=0$, $c=3$, take parts of size 1 indexed by $\emptyset$, the five pairs $\{i,i+1\}$ with indices modulo 5, and the sets $[5]\setminus\{0\}$ and $[5]\setminus\{2\}$, and make all other parts empty. The five pairs induce a $5$-cycle. The empty set is non-adjacent to the two 4-sets, whose only common neighbour among the five pairs is $\{3,4\}$. Thus this 8-vertex expansion contains no $K_4$ and  $\K{\barL}{4}\ge 8$. Conversely, every expansion of $\barL$ is $\OK_3$-free and if we require that it has no $K_4$ then its order must be strictly smaller than the Ramsey number $R(4,3)=9$. Therefore $\K{\barL}{4}=8$.

Let us show the opposite direction in the remaining cases, writing $t:=5m+c+1$. The case $t=1$ follows from the definition, so assume that $t\ge 2$. Let $G$ be an arbitrary expansion of $\barL$ with  $k:=16m+b+1$ vertices. We have to show that $G$ has a clique of order $t$. Suppose that this is false. 
Lemma~\ref{lm:General} (with $\omega=d=5$) implies that $16(t-1)\ge 5k$, that is, $16c\ge 5(b+1)$. This gives a contradiction for $c\in \{0,3,4\}$.

The following argument will be used to prove the lemma in the remaining cases $c=1,2$. Let $V_X$ be a part of the maximum size, denote it by $m+\mu$. Since $b\ge 2$, we have $\mu \ge1$. By the vertex transitivity of $\barL$,
we can assume that $X=\emptyset$. The set  of  non-neighbours of $X$ in $\barL$ is $\binom{[5]}4$, that is, it consists of all five 4-subsets of~$[5]$. It forms a 5-clique; thus the corresponding parts of the expansion have  at most $5m+c$ vertices in total. The set $\binom{[5]}2\subseteq V(\barL)$ of the neighbours of $X$  spans five copies of $K_4$ (each consisting of four pairs sharing a common element of $[5]$); also, each element of $\binom{[5]}2$ is in exactly two of these $4$-cliques. 
Each such $4$-clique extends to a $5$-clique by adding $\emptyset\in V(\barL)$, so the corresponding four parts contain at most $4m+c-\mu$ vertices. By summing this over all five such $4$-cliques, we get that the parts with indices in $\binom{[5]}2$ have at most $\frac52(4m+c-\mu)$ vertices in total. Thus, by $\mu\ge 1$, we have
 \[
 16m+b+1\le (m+\mu)+(5m+c)+\frac52(4m+c-\mu)=16m+ \frac{7c-3\mu}2\le 16m+\frac{7c-3}2.
 \]
 Thus $2b+2\le 7c-3$.  This,  for $c=1,2$, contradicts the definition of~$b$.\epf 
 
Note that the function $\K{\barL}{t}$ of Lemma~\ref{lm:L} is $16t/5+O(1)$ while $\K{\OK_{3}}{t}=3(t-1)$ by Lemma~\ref{lm:OK}.
Thus the combination of Theorems~\ref{th:ER}\ref{it:L} and~\ref{th:NewER} with  Proposition~\ref{pr:Reduction2} (for $H=\barL$ and $\C F=\{\OK_3\}$) applies only to finitely many pairs $(k,t)$. More specifically, for $t=5m+c+1$ with $m\ge 1$, we have $\K{\barL}{t}-3(t-1)=m+b(c)-3c\ge 0$. For $t=3,4,5$, the values of $\K{\barL}{t}$ are respectively $5,8,12$. Thus, for $t\ge 3$, precisely the pairs $(6,3)$ and $(9,4)$ satisfy~\eqref{eq:Reduction2}.
So we have the following corollary.

\begin{corollary}\label{cr:L} For $s\in\{6,7,8\}$ and $(k,t)\in\{(6,3),(9,4)\}$, it holds that $\er{s}{\BF{k,t}}$ is equal to the expression in~\eqref{eq:KsL} and every almost $\ER{s}{n,\BF{k,t}}$-extremal graph
is $o(n^2)$-close in the edit distance to a uniform expansion of $\barL$.\qed
\end{corollary}

Next, let us state the consequence of Theorem~\ref{th:ForbC5}. Combining it with Proposition~\ref{pr:Reduction2} (and using that $\min\{\K{C_5}{t},\K{\OK_3}{t}\}=\floor{5(t-1)/2}$ by Lemma~\ref{lm:OK} and \ref{lm:C5}), we obtain by~\eqref{eq:Reduction2} the following.

\begin{corollary}\label{cr:2} If $s\in\{4,5,6\}$ and integers $k\ge t\ge 2$ satisfy $2(t-1)< k\le 5(t-1)/2$
then  $\er{s}{\BF{k,t}}=2^{1-s}$ and every almost $\ER{s}{n,\BF{k,t}}$-extremal  graph  is $o(n^2)$-close in the edit distance to a uniform expansion of~$\OK_{2}$.\qed\end{corollary}

Let us turn to the consequences of Theorem~\ref{th:734}. Using the easy equality $\K{F\sqcup H}{t}=\K{F}{t}+\K{H}{t}$, where $F\sqcup H$ denotes the union of vertex-disjoint copies of $F$ and $H$, we have by Lemmas~\ref{lm:OK} and~\ref{lm:C5} that $\K{C_5'}{t}=\floor{5(t-1)/2}+t-1=\floor{7(t-1)/2}$. Thus we obtain the following corollary of Theorem~\ref{th:734}.

\begin{corollary}\label{cr:734} 
For any $k\ge t\ge 3$ with $3(t-1)< k\le 7(t-1)/2$, it holds that $\er{4}{\BF{k,t}}=1/27$ and every almost $\ER{4}{n,\BF{k,t}}$-extremal graph is $o(n^2)$-close in the edit distance to a uniform expansion of~$\OK_3$.\qed
\end{corollary}

Let us summarise what is known about the values of $\er{s}{\BF{k,t}}$ for $3\le t\le k\le 7$ and $s\ge 3$. Note that~\eqref{eq:1} gives that $\er{s}{\BF{k,t}}=1$ if $k\le 4$ (and any $t\ge 3$), or $k\in \{5,6\}$  and $t\ge 4$, or $k=7$ and $t\ge 5$. Table~\ref{ta:BB} lists, for the remaining solved triples of parameters with $s\ge 3$, the value of $\bb{s}{k,t}$ and the graph $H$ whose uniform expansions give the value of $\er{s}{\BF{k,t}}$. 
\begin{table}
\begin{center}
\begin{tabular}{c|c|c|c|c}
\hline
 $(k,t)$ &  
 $s$ & $\bb{s}{k,t}$ & $H$ & Reference \\
 \hline
 $(5,3)$ & $3,4,5,6$ & $2^{1-s}$ & $\OK_2$ & Corollaries~\ref{cr:1} and~\ref{cr:2}\\
 $(6,3)$ & $3$ & $2^{1-s}$ & $\OK_2$ & Corollary~\ref{cr:1} \\
 $(6,3)$ & $4,5$ & $(2^s-1)/5^{s-1}$ & $C_5$ & Corollary~\ref{cr:C5}\\
 $(6,3)$ & $6,7,8$ & the expression in~\eqref{eq:KsL} & $\barL$ & Corollary~\ref{cr:L}\\
 $(7,3)$ & $3,4$ & $3^{1-s}$ & $\OK_3$ & Corollaries~\ref{cr:1} and~\ref{cr:734}\\
 $(7,4)$ & $3,4,5,6$ & $2^{1-s}$ & $\OK_2$ & Corollaries~\ref{cr:1} and~\ref{cr:2}
\end{tabular}
\caption{Asymptotically solved cases of the Baumann--Briggs problem for $3\le t\le k\le 7$ and $s\ge 3$ with  $\er{s}{\BF{k,t}}<1$.}
\label{ta:BB}
\end{center}
\end{table}

All previous corollaries on the Baumann--Briggs problem give an asymptotic answer. Let us give one sample result describing the exact structure of large extremal graphs.

\begin{theorem}\label{th:53}
For all integers $s\in\{3,4,5,6\}$ and $t\ge 3$, if $n$ is sufficiently large then every $\ER{s}{n,\BF{2t-1,t}}$-extremal graph has a uniform expansion of $\OK_2$ as a spanning subgraph and it holds that $\ER{s}{n,\BF{2t-1,t}}=\binom{\floor{n/2}}s+\binom{\ceil{n/2}}s$.
\end{theorem}

\section{Proofs}

In order to prove Proposition~\ref{pr:Reduction2}, we will need the following auxiliary result. 

\begin{proposition}\label{pr:Reduction} For all integers $k\ge t\ge 2$, it holds as $n\to\infty$ that every $\BF{k,t}$-free $n$-vertex graph $G$ can be made $\BF{k,t}'$-free by changing at most $o(n^2)$ adjacencies.\end{proposition}
\bpf 
Let $\BF{k,t}''$ consist of those graphs in $\BF{k,t}'$ that have at most $k$ vertices. Every member of $\BF{k,t}'$ has an induced subgraph in $\BF{k,t}''$, obtained by retaining the vertices corresponding to the nonempty parts of a witnessing $k$-vertex expansion. Thus a graph is $\BF{k,t}'$-free if and only if it is $\BF{k,t}''$-free. 
Suppose that (some sequence of graphs) $G$ contradicts the statement. 
By the Induced Removal Lemma of Alon, Fischer, Krivelevich and Szegedy~\cite{AlonFischerKrivelevichSzegedy00}, there is $F\in\BF{k,t}''$, say with $V(F)=[m]$, such that $\pp(F,G)=\Omega(n^{m})$, where $n$ is the number of vertices in~$G$. 
Since $F\in \BF{k,t}'$, we can find
$H=\blow{F}{V_0,\dots,V_{m-1}}$, a $k$-vertex expansion of $F$ with no $t$-clique. Using, for example, the $m$-partite $m$-hypergraph result of Erd\H os~\cite{Erdos64}, we conclude that $G$ contains disjoint sets $U_0,\dots,U_{m-1}$ such that $|U_i|=|V_i|$ for $i\in [m]$ and every map $f:[m]\to V(G)$ with $f(i)\in U_i$ for $i\in [m]$ is an isomorphism of $F$ on its image. In other words, we have found a blowup of $F$ in $G$ with the same part sizes as in $H$ (where we do not stipulate anything about the edges inside a part $U_i$). Since $G$ is $\BF{k,t}$-free, there is a $t$-set $T\subseteq U_0\cup\dots\cup U_{m-1}$ that spans a clique in~$G$. But then if we choose a $t$-set $T'$ in $V_0\cup\dots\cup V_{m-1}$ with $|T'\cap V_i|=|T\cap U_i|$ for each $i\in[m]$ then it necessarily spans a $t$-clique in $H=\blow{F}{V_0,\dots,V_{m-1}}$ (because each $V_i$ induces a complete graph in~$H$). This contradicts our choice of~$H$.\epf

 \bpf[Proof of Proposition~\ref{pr:Reduction2}.] If we take an $n$-vertex expansion of $H$ with $\ER{s}{n,\C F}+o(n^s)$ copies of $K_s$, then the same expansion is admissible for the $\ER{s}{n,\BF{k,t}}$-problem since we assumed that $H\not\in\BF{k,t}'$. Thus $\ER{s}{n,\BF{k,t}}\le \ER{s}{n,\C F}+o(n^s)$.
 
In the other direction, each almost  $\ER{s}{n,\BF{k,t}}$-extremal graph $G$ can be changed into a $\BF{k,t}'$-free graph $G'$ by doing $o(n^2)$ adjacency edits by Proposition~\ref{pr:Reduction}. Of course, this affects the count of $K_s$-subgraphs by $o(n^s)$. Since $\C F\subseteq  \BF{k,t}'$, we have that
 \[
 \ER{s}{n,\C F}\le \ER{s}{n,\BF{k,t}'}\le \pp(K_s,G')\le \pp(K_s,G)+o(n^s)=\ER{s}{n,\BF{k,t}}+o(n^s).
 \] 
 
 Thus these two functions are within additive $o(n^s)$ from each other, giving~\eqref{eq:Reduction2=}. Also, the statement about almost extremal graphs clearly follows.\epf

\bpf[Proof of Theorem~\ref{th:NewER}.] 
The lower bound on $\er{8}{\OK_3}$ is proved via the standard application of the flag algebra method of Razborov~\cites{Razborov07,Razborov10}, using the \texttt{FlagAlgebraToolbox} package by the first author~\cite{Bodnar26fat} with commit \texttt{22b8765}. This method is well-known and is described in detail in e.g.~\cites{BaberTalbot11,SFS16,GilboaGlebovHefetzLinialMorgenstein22,JeongParkYang}.

In this particular case, we had to use base flags with $8$ vertices, when there are $410$ different  triangle-free graphs up to isomorphism.

The obtained exact-arithmetic certificate also satisfies (after complementing all graphs so that we work with blowups instead of expansions) the sufficient condition of the second author, Slia\v can and Tyros~\PST{Theorem 7.1(i)} for the so-called
\emph{perfect $\OO R_{3,3,3}$-stability} with the uniform vector $(1/16,\dots,1/16)$ as the unique optimal vector of part ratios. We do not define this property here but observe that it implies all the remaining claims of the theorem and refer the reader to~\PST{} for further details. 
In order to apply the condition from~\PST{}, one has to specify a graph $\tau$ as an extra input (for which, in the complementary setting, we used $\tau=C_5'$, the $5$-cycle plus an isolated vertex) and check various assumptions which were formulated in such a way in~\PST{} that they can be verified by a computer. This verification is done by our scripts. In fact, some of these properties can be verified without using computer and this was already done in~\PV{Claim~4.9}, which states that there is a unique way (up to an automorphism of $\barL$) to embed the complement of $C_5'$ into an expansion of $\barL$ and each part of the expansion is uniquely identified by its attachment to the image of $C_5'$.  

The Jupyter notebook generating the certificates and verifying all stated properties (including the criterion for perfect stability from \PST{}) and the certificates  are included with the arXiv submission as ancillary files and are also available at \url{https://github.com/bodnalev/supplementary_files/tree/main/erdos_baumann_briggs} (with all certificates available in the \texttt{certificates} folder). 
For the floating-point calculation, the high-precision SDP solver SDPA-QD was used \cite{Nakata10}. 
The generated certificate for $\er{8}{\OK_3}$ is stored as \texttt{er83\_cert.pickle}.\epf

\bpf[Proof of Theorems~\ref{th:ForbC5} and~\ref{th:734}.] Recall that we have to prove various claims about  the $\ER{s}{n,\C F}$-problem where  $\C F=\{\OK_{\ell+1},C_5\}$, $\ell=2$ and $s\in\{4,5,6\}$ in Theorem~\ref{th:ForbC5} and $\C F=\{\OK_{\ell+1},C_5'\}$, $\ell=3$ and  $s=4$  in Theorem~\ref{th:734}.

Again, we use flag algebras to prove that $\er{s}{\C F}$ is at least (and thus equal) to the value coming from uniform expansions of $\OK_{\ell}$ in each stated case. The corresponding certificates were generated using the \texttt{FlagAlgebraToolbox} package~\cite{Bodnar26fat}. All SDPs were solved using the high-precision solver SDPA-QD.
For Theorem~\ref{th:ForbC5} we had to use base flags with $N:=s+2$ vertices (namely $N=6,7,8$ for $s=4,5,6$ respectively). The rounded certificates are stored as
\texttt{BB\_e3c5\_sX\_cert.pickle}, where \texttt{X} is the value of $s$.
For Theorem~\ref{th:734} we use base flags with $N:=7$ vertices and the rounded certificate is stored as \texttt{BB\_e4c5k1\_s4\_cert.pickle}.

Since the arguments proving the remaining claims of both theorems can be made rather similar, we present a single proof.

Let us prove the claimed stability property. Let $n\to\infty$ and take any almost $\ER{s}{n,\C F}$-extremal graph $G=(V,E)$. We have to show that $G$ is $o(n^2)$-close in the edit distance to a uniform expansion of~$\OK_\ell$.
Each certificate in fact proves for some constant $c>0$ that the density of $K_s$ minus $c$ times the induced density of $K_N^-$ (the $N$-clique with one edge missing) is at least $\ell^{1-s}+o(1)$. It follows that $\pp(K_N^-,G)=o(n^N)$. 
By the Induced Removal Lemma~\cite{AlonFischerKrivelevichSzegedy00}, applied to the finite family $\C F\cup\{K_N^-\}$, we can change $o(n^2)$ adjacencies in $G$ so that it has no induced copy of $K_N^-$ and is still $\C F$-free. Let $S\subseteq V$ be a set spanning a maximum clique in $G$. Since $G$ is $\OK_4$-free, Ramsey's Theorem implies that $|S|$ goes to infinity. By the maximality of $S$, every vertex $u\in V\setminus S$ has at least one non-neighbour in~$S$. Since $G$ is $K_N^-$-free,  $u$ has at most $N-2$ neighbours in $S$. 


Suppose first that $\ell=2$ (that is, we prove Theorem~\ref{th:ForbC5}). Since $|S|>2(N-2)$, every two vertices outside of $S$ have a common non-neighbour in $S$. As  $G$ is $\OK_3$-free, we conclude that $V\setminus S$ spans a clique. Thus $G$ is a union of two cliques plus at most $(N-2)\,|V\setminus S|=o(n^2)$ edges across. Since the $K_s$-density in $G$ is $2^{1-s}+o(1)$, an easy calculation (e.g.\ using strict convexity of $x\mapsto x^s$) shows each clique must have $(1/2+o(1))n$ vertices, giving the required stability property.

Now, let $\ell=3$  (that is, we prove Theorem~\ref{th:734}). Since $|S|>3(N-2)$, every three vertices outside of $S$ have a common non-neighbour in $S$ and, by the $\OK_4$-freeness of $G$, thus span at least one edge. Thus $G':=G-S$ is $\OK_3$-free. Let $S'\subseteq V\setminus S$ span a maximum clique in~$G'$. We have that $|V\setminus S|=\Omega(n)$ (as otherwise the $K_s$-density in $G$ would be $1-o(1)$) and that $|S'|\to\infty$ by Ramsey's Theorem. 
As before, we argue that every vertex in $V'':= V\setminus(S\cup S')$ sends at most $N-2$ edges to $S'$. Therefore, if we take any two vertices $u,v\in V''$ then their common non-neighbourhoods in each of $S$ and $S'$ are at least, say, $N-1$ and $1$ respectively, and thus there is at least one non-edge between these two sets. Since $G$ is $\OK_4$-free, we conclude that $V''$ spans a clique in~$G$. Thus $G$ is a union of three cliques plus at most $2(N-2)n=o(n^2)$ edges across. Since the $K_s$-density in $G$ is $3^{1-s}+o(1)$ (for $s=4$), it follows that each clique has $(1/3+o(1))n$ vertices, as desired.

It remains to analyse the structure of large extremal graphs. Let $n$ be sufficiently large and let $G=(V,E)$ be an arbitrary graph which is $\ER{s}{n,\C F}$-extremal, that is, $G$ is
an $\C F$-free $n$-vertex graph with $\pp(K_s,G)=\ER{s}{n,\C F}$. Let $V_0\cup\dots\cup V_{\ell-1}$ be a partition of $V$ into $\ell$ parts that minimises $|W|$, where $W$ consists of non-edges inside a part and edges across the parts; we call the pairs in $W$ \emph{wrong}. 

By stability, there is a partition with $o(n^2)$ wrong pairs, so our choice gives $|W|=o(n^2)$. Also, the total number of edges in $G$ is $(1/\ell+o(1))\binom n2$. Therefore
\[
\sum_{i=0}^{\ell-1}\binom{|V_i|}{2}=|E(G)|+o(n^2)=(1/\ell+o(1))\binom n2.
\]
Since $\sum_{i=0}^{\ell-1}|V_i|=n$, it follows that $\sum_{i=0}^{\ell-1}(|V_i|-n/\ell)^2=o(n^2)$, so each part has $(1/\ell+o(1))n$ vertices.

The average number of $K_s$-subgraphs per vertex is $(s/n)\cdot \pp(K_s,G)= (\ell^{1-s}+o(1))\binom{n-1}{s-1}$. Also, the numbers of $K_s$-subgraphs per any two vertices $u$ and $v$ differ at most by $\binom{n-2}{s-2}$, since we can delete $u$ and add a new vertex $u'$ adjacent to $v$ whose neighbourhood in $V\setminus\{u,v\}$ is the same as that of $v$; this operation keeps the graph $\C F$-free, since none of the forbidden graphs has two adjacent vertices with the same neighbours outside that pair, and thus cannot decrease the number of copies of~$K_s$.
Thus each vertex is in $(\ell^{1-s}+o(1))\binom{n-1}{s-1}$ $s$-cliques.

Take any vertex $v\in V_i$ of $G$. For $j\in [\ell]$, write $d(v,V_j):=|N_G(v)\cap V_j|$. By the minimality of $|W|$, moving $v\in V_i$ to $V_j$ cannot improve the partition, so $2d(v,V_i)-|V_i|+1\ge 2d(v,V_j)-|V_j|$. Since the part sizes differ by $o(n)$, it follows that $d(v,V_i)\ge d(v,V_j)+o(n)$ for every $j$, uniformly in $v$. We claim that $v$ is adjacent to $(1/\ell+o(1))n$ vertices of $V_i$ (that is, to almost every vertex in its part). Indeed, otherwise $v$ would have $\Omega(n)$ non-neighbours in every part (since the part sizes differ by $o(n)$ only)
and the corresponding $\ell$ sets of non-neighbours of $v$ would contain a copy of $\OK_\ell$ (since $G$ has $o(n^2)$ cross edges), which together with $v$ would give a copy of $\OK_{\ell+1}$ in $G$, a contradiction. It follows that $v$ is in $(\ell^{1-s}+o(1))\binom{n-1}{s-1}$ $s$-cliques that lie entirely inside its part $V_i$ and, therefore, is in $o(n^{s-1})$ $s$-cliques that intersect some other part. The last property implies that $v$ has $o(n)$ neighbours in every part different from~$V_i$.

Since $v$ was arbitrary, we see that $W$, when viewed as a graph, has maximum degree $\Delta(W)=o(n)$. 

Let us show that each part $V_i$ spans a clique in $G$. Indeed, if we had some non-adjacent distinct vertices $u_0$ and $u_1$ in, say, $V_0$ then we could, iteratively for $i=1,\dots,\ell-1$, pick $u_{i+1}\in V_i$ non-adjacent to each of $u_0,\dots,u_i$ since at most $\ell\Delta(W)=o(n)$ vertices of $V_i$ are excluded; this would produce a copy of $\OK_{\ell+1}$ in $G$, a contradiction.

We have just shown that $G$ has the expansion $G':=\blow{\OK_{\ell}}{V_0,\dots,V_{\ell-1}}$ as a spanning subgraph. Since $G'$ is $\C F$-free (and $G'\subseteq G$), it is also $\ER{s}{n,\C F}$-extremal. It easily follows that the part sizes differ at most by one (as otherwise moving a vertex from a largest part to a smallest one would strictly decrease the number of $s$-cliques). This finishes the proof of Theorems~\ref{th:ForbC5} and~\ref{th:734}.\epf

\bpf[Proof of Theorem~\ref{th:53}.] 
The argument starts similarly to the proof of Theorem~\ref{th:ForbC5}.

Take an extremal graph $G$ of sufficiently large order $n$ and choose a partition $V(G)=V_0\cup V_1$ minimising $|W|$, where $W:=E(G)\bigtriangleup E(\blow{\OK_2}{V_0,V_1})$. Stability, which follows from Corollary~\ref{cr:1} for $s=3$ and Corollary~\ref{cr:2} for $s\in\{4,5,6\}$, gives $|W|=o(n^2)$. As in the preceding proof, each part has size $(1/2+o(1))n$, and moving a vertex between the parts shows that $d(v,V_i)\ge d(v,V_{1-i})+o(n)$ for $v\in V_i$, uniformly in~$v$.

Suppose that a vertex $v$ of $G$ has $\Omega(n)$ non-neighbours in its part $V_i$. Then $v$ also has $\Omega(n)$ non-neighbours in the other part $V_{1-i}$. Since $G$ has $o(n^2)$ edges between $V_0$ and $V_1$, we can find $(t-1)$-sets $X_i\subseteq V_i$ of non-neighbours of $v$  for $i=0,1$ with $G$ having no edges between $X_0$ and~$X_1$. But then the $(2t-1)$-set $X_0\cup X_1\cup\{v\}$ spans no $t$-clique, a contradiction. 

Thus the maximum number of non-neighbours of a vertex inside its own part is $o(n)$.

Let $U$ consist of the vertices having at least $n/8$ neighbours in the other part. We claim that every missing edge inside a part has an endpoint in $U$. Suppose that this is false and take non-adjacent $u,v\in V_i\setminus U$. Since there are $o(n^2)$ cross edges, all but $o(n)$ vertices have $o(n)$ neighbours in the other part. Extend $\{u,v\}$ to a $t$-set $X_i\subseteq V_i$ using $t-2$ such vertices. At least $|V_{1-i}|-n/4+o(n)=\Omega(n)$ vertices in the other part send no edges to $X_i$. Choosing a $(t-1)$-set $X_{1-i}$ among them gives a $(2t-1)$-set $X_0\cup X_1$ spanning no $t$-clique, a contradiction.

Let $M$ be the number of missing edges inside the parts. By the above, $M\le |U|o(n)$. Each vertex of $U$ lies in $\Omega(n^{s-1})$ $s$-cliques whose other vertices lie in the opposite part, since its at least $n/8$ neighbours there span all but $o(n^2)$ of their possible edges. Each such clique is counted at most $s$ times, so there are at least $c|U|n^{s-1}$ cross-part $s$-cliques for some fixed $c>0$. On the other hand, the missing internal edges destroy at most $M\binom{n-2}{s-2}=o(|U|n^{s-1})$ of the $s$-cliques in $G':=\blow{\OK_2}{V_0,V_1}$. If $U\ne\emptyset$, this gives $\pp(K_s,G)>\pp(K_s,G')$ for all sufficiently large $n$, contradicting extremality because $G'$ is $\BF{2t-1,t}$-free. Thus $U=\emptyset$, and the claim gives $M=0$.

Thus $G$ has $G'$ as a spanning subgraph. Since $G'$ is also extremal, its part sizes differ by at most one: otherwise moving a vertex from the larger part to the smaller one would strictly decrease the number of $s$-cliques.\epf

Let us remark that  all extremal graphs of large order $n$ in Theorems~\ref{th:ForbC5}, \ref{th:734} and \ref{th:53} can be described as follows: for the appropriate $\ell$ and $s$, we take a uniform $n$-vertex expansion of $\OK_\ell$ and add an arbitrary set of edges so that no new $s$-cliques are created.

\section{Concluding remarks}

We did not try to determine the exact value of the Baumann--Briggs function apart from Theorem~\ref{th:53} (when $k=2t-1$), since the proofs and calculations for general $(k,t)$ will probably be rather messy. For example, consider the simplest case of Corollary~\ref{cr:1} when $\OK_2\not\in \BF{k,t}'$ and $\OK_3\in\BF{k,t}'$ (when we know that uniform expansions of $\OK_2$ are almost extremal for $s=3$). Here, already for $k=2t\ge 6$, we have the freedom of removing edges (namely, a star) from one part of a $\OK_2$-expansion without creating any forbidden subgraph. The situation when  $k-2t$ assumes values larger than 0 is more complicated since we have an option of removing edges from both parts.
Dealing with these complications will probably require technical work rather than new insights, so we do not pursue this direction here.

Note that the conclusion of Theorem~\ref{th:734} will not be true if we forbid $\OK_4$ only, since (non-uniform) 
expansions of the unique 8-vertex $(3,4)$-Ramsey graph with 10 edges, denoted by $R_{3,4}$, show that $\er{4}{\OK_4}\le ({14\cdot 2^{1/3}-11})/{192}$ (see \PV{Pages 932--933}), which is strictly less than $1/27$. Thus e.g.\ the case $(k,t,s)=(7,3,4)$ of the Baumann--Briggs problem is an example where the determination of $\bb{s}{k,t}$ does not reduce to some instance of the Erd\H os $\er{s}{\OK_\ell}$-problem via Proposition~\ref{pr:Reduction2}.

In an earlier version of\hide{this} the arXiv pre-print, we asked if the conclusion of Theorem~\ref{th:ForbC5} holds for every $s\ge 7$. As it was pointed out to us by Sergey Norin, this is not the case: uniform expansions of, for example, $\O{C}_7$ show that $\er{s}{\{\OK_3,C_5\}}\le 7\cdot (3/7)^s - 7\cdot (2/7)^s$ which is smaller than $2^{1-s}$ for $s \geq 8$. 

\section*{Acknowledgements}

The authors were supported 
by ERC Advanced Grant 101020255. The authors would like to thank the anonymous referee for pointing out errors in the initial version of the paper.

\hide{
\section*{Macros}

$\ER{s}{n,\OK_m}$ \verb$\ER{s}{n,\OK_m}$ : Erdos fn

$\er{s}{\OK_m}$ \verb$\er{s}{\OK_m}$ : scaled Erdos fn


$\bb{s}{k,t}$ \verb$\bb{s}{k,t}$ : scaled Baumann--Briggs fn

$\BF{k,t}$ \verb$\BF{k,t}$ : Baumann--Briggs family

$\BF{k,t}'$ \verb$\BF{k,t}'$ : $H$ st some $k$-vertex expansion has no $t$-clique


$\blow{H}{V_0,\dots,V_{m-1}}$ \verb$\blow{H}{V_0,\dots,V_{m-1}}$ : expansion

$\OK_m$ \verb$\OK_m$ : empty $m$-vertex graph

$\R$ \verb$\R$ : 13-vertex 35-Ramsey graph

$\barL$ \verb$\barL$ : complement of Clebsch
}

\bibliography{bibexport}
\end{document}